\documentclass{IEEEtran}
\usepackage{amsmath,amssymb,amsfonts,amsthm}
\usepackage{graphicx}
\usepackage{textcomp}
\usepackage{xcolor}
\def\BibTeX{{\rm B\kern-.05em{\sc i\kern-.025em b}\kern-.08em
    T\kern-.1667em\lower.7ex\hbox{E}\kern-.125emX}}

\usepackage{bm}
\usepackage[caption=false, font=footnotesize]{subfig}
\usepackage{url}

\newtheorem{lemma}{Lemma}
\newtheorem{theorem}{Theorem}
\newtheorem{remark}{Remark}

\begin{document}
\title{Stability and Convergence of a Randomized Model Predictive Control Strategy}
\author{Dani\"el W. M. Veldman, Alexandra Borkowski, and Enrique Zuazua 
\thanks{Submitted to the editors on \today. E. Zuazua has been funded by the Alexander von Humboldt-Professorship program, the ModConFlex Marie Curie Action, HORIZON-MSCA-2021-DN-01, the COST Action MAT-DYN-NET, the Transregio 154 Project ``Mathematical Modelling, Simulation and Optimization Using the Example of Gas Networks'' of the DFG, grants PID2020-112617GB-C22 and TED2021-131390B-I00 of MINECO (Spain), and by the Madrid Goverment -- UAM Agreement for the Excellence of the University Research Staff in the context of the V PRICIT (Regional Programme of Research and Technological Innovation).}
\thanks{D.~W.~M. Veldman and E. Zuazua are with the Chair for Dynamics, Control, Machine Learning and Numerics -- Alexander von Humboldt Professorship, Department of Mathematics, Friedrich-Alexander Universit\"at Erlangen-N\"urnberg, 91058, Erlangen, Bavaria, Germany (e-mail: daniel.wm.veldman@fau.de, enrique.zuazua@fau.de). }
\thanks{A. Borkowski is with the Department of Mathematics, King's College London, Strand, London, WC2R 2LS, United Kingdom (e-mail: alexandra.borkowski@kcl.ac.uk).}
\thanks{E. Zuazua is also with the Chair in Computational Mathematics, Fundaci\'on Deusto, Av. de las Universidades 24, 48007, Bilbao, Spain and Departamento de Matem\'aticas, Universidad Autónoma de Madrid, 28049, Madrid, Spain.
}
}

\maketitle

\begin{abstract}
RBM-MPC is a computationally efficient variant of Model Predictive Control (MPC) in which the Random Batch Method (RBM) is used to speed up the finite-horizon optimal control problems at each iteration. In this paper, stability and convergence estimates are derived for RBM-MPC of unconstrained linear systems. 
The obtained estimates are validated in a numerical example that also shows a clear computational advantage of RBM-MPC. 
\end{abstract}


\section{Introduction}
Model Predictve Control (MPC) is a well-established and widely used method to control complex dynamical systems, see, e.g., \cite{mayne2000, grune2017, rawlings2019} for an overview of the large body of research in this area. MPC requires the real-time solution of a sequence of optimal control problems (OCPs) on a finite time horizon, which can be computationally demanding. This is for example the case when the model is the result of the (spatial) discretization of a Partial Differential Equation (PDE) or in the simulation of interaction particle systems. 

One recently-proposed numerically-efficient approximation method is the Random Batch Method (RBM) \cite{jin2020}, which is closely related to the stochastic algorithms like Stochastic Gradient Descent (SGD). In the RBM-dynamics, 
a random subset/batch of 
interconnections between Degrees of Freedom (DOFs)
is considered during small time intervals.
This can reduce the computational cost significantly and leads to a good approximation of the original dynamics when these time intervals are chosen sufficiently small, see, e.g., \cite{jin2020}. Recently, this idea has been extended to infinite-dimensional systems
\cite{eisenmann2022}. 
RBM-constrained OCPs have been analyzed in \cite{rbm}.

The RBM can be used to speed up the solution of the finite-horizon OCPs in MPC. The feedback nature of MPC also creates more robustness against the accumulating error in the RBM approximation (see, e.g., \cite{rbm}). The effectiveness of this combination of MPC with RBM (RBM-MPC) for nonlinear interacting particle systems has been demonstrated in \cite{ko2021model}, but a rigorous stability and convergence analysis is still missing. 

The RBM in RBM-MPC fulfills a similar role as the Reduced Order Models (ROM) in MPC based on ROMs. 
There has been research on stability guarantees  for MPC based on ROMs in constrained linear systems, see, e.g., \cite{loehning2014, lorenzetti2022}. The RBM is typically easier to apply than ROM techniques, but the analysis of RBM-MPC is nontheless involved due to the stochasticity introduced by the RBM.%

In this paper, we provide the first rigorous analysis of the RBM-MPC algorithm. Our analysis is limited to the unconstrained linear quadratic setting and thus extends the open-loop analysis from \cite{rbm} to a closed-loop setting. The obtained error estimates demonstrate the influence of the different parameters in RBM-MPC on the expected performance, and the obtained convergence rates are validated in a numerical example. 

The remainder of this paper is structured as follows. The RBM-MPC algorithm is presented in Section \ref{sec:rbm_mpc}. After the introduction of preliminary estimates and notation in Section \ref{sec:preliminaries}, the stability and convergence of RBM-MPC are proven in Sections \ref{sec:stab} and \ref{sec:conv}, respectively. The convergence rates are validated in a numerical example in Section \ref{sec:numerical}. Finally, conclusions and perspectives are presented in Section \ref{sec:conclusions}. 

We will use the following notation. 
The (Euclidean) norm of a vector $\bm{x} \in \mathbb{R}^n$ is $|\bm{x}| = \sqrt{\bm{x}^\top\bm{x}}$. For a matrix $M \in \mathbb{R}^{n \times m}$, $\| M \| = \sup_{|\bm{x}| = 1} |M\bm{x}|$. 
For symmetric $M \in \mathbb{R}^{n \times n}$, $M \succcurlyeq 0$ or $M \succ 0$ indicates that $M$ is positive semi-definite or positive definite, respectively. For $M \succcurlyeq 0$, $|\bm{x}|_M = \sqrt{\bm{x}^\top M \bm{x}}$. 

\section{The RBM-MPC algorithm} \label{sec:rbm_mpc}


The RBM-MPC algorithm analyzed in this paper is a way to approximate the control $\bm{u}^*_\infty(t)$ that minimizes
\begin{equation}
    J_\infty(\bm{u}) = \int_0^\infty \left(|\bm{x}(t)|^2_Q + |\bm{u}(t)|_W^2 \right) \ \mathrm{d}t, \label{eq:Jinf}
\end{equation}
subject to the dynamics
\begin{equation}
    \dot{\bm{x}}(t) = A \bm{x}(t) + B \bm{u}(t), \qquad \qquad \bm{x}(0) = \bm{x}_0,  \label{eq:dyn_x}
\end{equation}
where the state $\bm{x}(t)$ evolves in $\mathbb{R}^n$ starting from the initial condition $\bm{x}_0$, the control $\bm{u}(t)$ evolves in $\mathbb{R}^m$, $0 \preccurlyeq Q \in \mathbb{R}^{n \times n}$, $0 \prec W \in \mathbb{R}^{m \times m}$, $A \in \mathbb{R}^{n\times n}$, and $B \in \mathbb{R}^{n\times m}$. 
It is assumed that $(A,B)$ is stabilizable and $(A,Q)$ is detectable.  

%


The OCP \eqref{eq:Jinf}--\eqref{eq:dyn_x} can be approximated using MPC. 
In MPC, two parameters arise: the prediction horizon $T$ and the shorter control horizon $\tau$. Set $\tau_i := i \tau$ (with $i \in \mathbb{N}$) and let $\bm{u}^*_{T}(t;\bm{x}_{i-1},\tau_{i-1})$ and $\bm{x}^*_{T}(t;\bm{x}_{i-1},\tau_{i-1})$ denote the control and state trajectory that minimize
\begin{multline}
     J_{T}(\bm{u}_{T}; \bm{x}_{i-1}, \tau_{i-1}) = |\bm{x}_{T}(\tau_{i-1}+T)|_F^2 \\ +\int_{\tau_{i-1}}^{\tau_{i-1} + T} \left( |\bm{x}_{T}(t)|_Q^2 + |\bm{u}_{T}(t)|^2_W \right) \ dt, \label{eq:JT}
\end{multline}
where $F \succcurlyeq 0$ and  $\bm{x}_{T}(t)$ fulfills for $t \in [\tau_{i-1}, \tau_{i-1} + T]$
\begin{align}
    \dot{\bm{x}}_{T}(t) = A\bm{x}_{T}(t) + B\bm{u}_{T}(t),\;\;\;\;\;\;\;\;\bm{x}_{T}(\tau_{i-1}) = \bm{x}_{i-1}. \label{eq:dyn_xT}
\end{align}
When $n$ is large and $A$ not sparse, finding $\bm{u}^*_{T}(t;\bm{x}_{i-1},\tau_{i-1})$ and $\bm{x}^*_{T}(t;\bm{x}_{i-1},\tau_{i-1})$ is computationally demanding. We therefore replace $A$ by a randomized sparser matrix $A_R$.  

The randomized matrix $A_R(\bm{\omega}_i,t)$ is constructed as follows. First, $A$ is written as the sum of \emph{sparse} submatrices $A_m$
\begin{align}
    A = \sum_{m=1}^M A_m. \label{eq:Asplit}
\end{align}

Next, the subsets of $\{1,2, \ldots, M \}$ are enumerated as $S_1, S_2, \dots, S_{2^M}$ and a probability $p_\omega \in [0,1]$ is assigned to each subset $S_\omega$ ($\omega \in \{1,2,\ldots, 2^M \}$) such that $\sum_{\omega} p_\omega = 1$. 

The time interval $[0, T]$ is divided into $K$ time intervals of equal length $h$. For each of the $K$ time intervals, an element $\omega_{i,k} \in \{1,2, \ldots, 2^M \}$ of the vector $\bm{\omega}_i$ is selected according to the probabilities $p_\omega$. The matrix $A_R$ is now defined as follows
\begin{align}
    A_R(\bm{\omega}_i,t) = \sum_{m \in S_{\omega_{i,k}}} \frac{A_m}{\pi_m}, \qquad
    t \in [(k-1)h,kh). 
    \label{eq:def_AR}
\end{align}
The scaling factors $\pi_m$ are defined such that the expected value of $A_R(\bm{\omega}_i,t)$ is equal to $A$. In particular, $\pi_m$ denotes the probability of having the index $m$ in the selected subset
\begin{equation}
    \pi_m := \sum_{\omega \in \{ \omega' \in \{ 1,2,\dots,2^M\} \mid m \in S_{\omega'} \} } p_\omega. \label{eq:pim}
\end{equation}
The definition of $A_R$ thus requires that the probabilities $p_\omega$ are selected such that $\pi_m > 0$ for all $m \in \{1,2, \ldots, M \}$. 

The dynamics generated by $A_R(\bm{\omega}_i,t)$ is in expectation close to the dynamics generated by $A$ for $h$ sufficiently small (see \cite{rbm}) and replacing $A$ by $A_R(\bm{\omega}_i,t)$ reduces the computational cost when $A_R(\bm{\omega}_i,t)$ is much sparser than $A$. 
Consider therefore the control $\bm{u}^*_{R}(\bm{\omega}_i, t;\bm{x}_{i-1},\tau_{i-1})$ and state trajectory $\bm{x}^*_{R}(\bm{\omega}_i, t;\bm{x}_{i-1},\tau_{i-1})$ that minimize
\begin{multline}
     J_R(\bm{u}_R; \bm{\omega}_i, \bm{x}_{i-1}, \tau_{i-1})  =  |\bm{x}_{R}(\bm{\omega}_i, \tau_{i-1}+T)|_F^2 \\ 
     +\int_{\tau_{i-1}}^{\tau_{i-1}+T} \left( |\bm{x}_{R}(\bm{\omega}_i,t)|^2_Q + |\bm{u}_{R}(t)|_W^2 \right) \ dt, 
     \label{eq:JR}
\end{multline}
where $\bm{x}_{R}(t)$ fulfills for $t \in [\tau_{i-1}, \tau_{i-1}+T]$
\begin{multline}
    \dot{\bm{x}}_{R}(\bm{\omega}_i, t) = A_R(\bm{\omega}_i, t - \tau_{i-1})\bm{x}_{R}(\bm{\omega}_i, t) + B\bm{u}_{R}(t),\\
    \bm{x}_{R}(\bm{\omega}_i,\tau_{i-1}) = \bm{x}_{i-1}. \label{eq:dyn_xR}
\end{multline}

It has been proven in \cite{rbm} that $\bm{u}^*_{R}(\bm{\omega}_i, t;\bm{x}_{i-1},\tau_{i-1})$ is (in expectation) close to $\bm{u}^*_T(t; \bm{x}_{i-1}, \tau_{i-1})$ for $h$ small enough, see also Section \ref{sec:preliminaries}. Because $\bm{u}_R^*(\bm{\omega}_i, t; \bm{x}_{i-1}, \tau_{i-1})$ is used to control the dynamics generated by $A$, consider also the solution $\bm{y}_R^*(\bm{\omega}_i,t; \bm{x}_{i-1}, \tau_{i-1})$ of
\begin{multline}
    \dot{\bm{y}}_R^*(\bm{\omega}_i, t) = A \bm{y}_R^*(\bm{\omega}_i, t) + B \bm{u}^*_R(\bm{\omega}_i, t; \bm{x}_{i-1}, \tau_{i-1}), \\ \bm{y}_R^*(\bm{\omega}_i, \tau_{i-1}) = \bm{x}_{i-1}. \label{eq:dyn_yR}
\end{multline}
where $\bm{y}_R^*(\bm{\omega}_i, t)$ denotes $\bm{y}_R^*(\bm{\omega}_i,t; \bm{x}_{i-1}, \tau_{i-1})$ for brevity. 

The RBM-MPC algorithm now computes the control $\bm{u}_{R-M}(t)$ and state trajectory $\bm{x}_{R-M}(t)$ on $[0,\infty)$ as follows. 
\begin{enumerate}
    \item Initialize $\bm{x}_{R-M}(0) = \bm{x}_0$ and $i = 1$.
    \item Select a random vector $\bm{\omega}_i \in \{1,2,\ldots,2^M \}^K$.  
    \item Compute $\bm{u}_R^{*}(\bm{\omega}_i, t;\bm{x}_{R-M}(\tau_{i-1}), \tau_{i-1})$ and \\ $\bm{y}_R^{*}(\bm{\omega}_i, t;\bm{x}_{R-M}(\tau_{i-1}), \tau_{i-1})$ on $[\tau_{i-1}, \tau_{i-1} + T]$.
    \item Set $\bm{u}_{R-M}(t) = \bm{u}_R^{*}(\bm{\omega}_i, t;\bm{x}_{R-M}(\tau_{i-1}), \tau_{i-1})$ and \\ 
    $\bm{x}_{R-M}(t) = \bm{y}^*_R(t,\bm{x}_{R-M}(\tau_{i-1}), \tau_{i-1})$ on $[\tau_{i-1}, \tau_{i}]$.
    \item Set $i = i + 1$ and got to Step 2. 
\end{enumerate}
Note that RBM-MPC reduces to standard MPC when $A_R(\bm{\omega}_i,t) = A$ and that $\bm{x}_{R-M}(\tau_i)$ depends on the previously selected sequences $\bm{\omega}_j$ with $j \leq i$, which are denoted by
\begin{equation}
    \Omega_{i} := (\bm{\omega}_1, \bm{\omega}_2, \ldots, \bm{\omega}_{i}) \in \{1,2,\ldots, 2^M \}^{iK}. \label{eq:def_Omegai}
\end{equation}

The construction of the matrix $A_R(\bm{\omega}_i,t)$ leaves freedom in the choice of the submatrices $A_m$, the probabilities $p_\omega$, and the grid spacing $h$.
As for the submatrices $A_m$, splittings of the form \eqref{eq:Asplit} are standard in operator-splitting methods, which are well-established in the numerical analysis, see, e.g., \cite{quarteroni1994}. The specific choice of the $A_m$'s is often guided by physical insight. In many finite-dimensional examples, each $A_m$ represents an interaction between two Degrees of Freedom (DOFs) so that $M \leq n(n-1)/2$. 
Regarding the grid spacing $h$, note that the estimates in Theorems \ref{thm:1} and \ref{thm:2} below are proportional to $\sqrt{h \mathrm{Var}[A_R]}$, where
\begin{align}
    \mathrm{Var}[A_R] &:= \sum_{\omega = 1}^{2^M} \left\| A - \sum_{m \in S_{\omega}} \frac{A_m}{\pi_m}  \right\|^2 p_\omega. \label{eq:def_varA}
\end{align}
Reducing $\mathrm{Var}[A_R]$ thus enables us to use a larger step size $h$. Finally, note that assigning nonzero probabilities $p_\omega$ to larger subsets $S_\omega$ reduces $\mathrm{Var}[A_R]$, but will also make $A_R(\bm{\omega}_i,t)$ less sparse and thus potentially increases the computational cost
, see \cite[Section 2.3]{rbm} for further discussions and examples.

Error estimates for the RBM, as in Section \ref{sec:preliminaries} and in Theorems \ref{thm:1} and \ref{thm:2}, require a uniform quasi-dissipativity bound on $A_R$ in the tradition of \cite{lumer1961}, i.e. we fix a $\mu_R \geq 0$ such that
\begin{equation}
\bm{x}^\top A_R(\bm{\omega}_i, t) \bm{x} \leq \mu_R |\bm{x}|^2, \label{eq:def_muR}
\end{equation}
for all $\bm{x} \in \mathbb{R}^n$, $\bm{\omega}_i \in \{1,2, \ldots, 2^M \}^K$, and $t \in [0, T]$. Note that this condition implies that the eigenvalues of the symmetric part of $A_R(\bm{\omega}_i, t)$ do not exceed $\mu_R$.

\begin{remark} \label{rem:dissipative}
Note that $\mu_R = 0$ when $\bm{x}^TA_m\bm{x}\leq 0$ for all $\bm{x} \in \mathbb{R}^n$ and all $m \in \{1, 2, \ldots, M\}$, i.e. when all $A_m$ are \emph{dissipative}. The latter condition can be achieved in many examples, see, e.g., Section \ref{sec:numerical} and \cite[Section 4]{rbm}. 
\end{remark}

\begin{remark} \label{PDE}
    Condition \eqref{eq:def_muR} readily extends to a setting in which the $A_m$ are quasi-dissipative unbounded operators on a Banach space, see, e.g.\ \cite{lumer1961}. However, the appearance of the operator norm $\| \cdot \|$ in \eqref{eq:def_varA} is an indication that extending the RBM to such setting is not trivial, see, e.g.\ \cite{eisenmann2022}. 
\end{remark}

\section{Preliminary estimates} \label{sec:preliminaries}

In the following, $C$ denotes a constant depending only on $A$, $B$, $Q$, $W$, and $F$. The notation $C_T$ indicates that the constant also depends on $T$. The constants $C$ and $C_T$ may vary from expression to expression, e.g.\ $(\| A \| + T)C_T \leq C_T$. Because we are interested in the limit $h \mathrm{Var}[A_R] \rightarrow 0$, we will only consider the lowest power of $h \mathrm{Var}[A_R]$ in our estimates. 

The following lemma now directly follows from \eqref{eq:def_muR}. 
\begin{lemma} \label{lem:bound_muRxR}
The solution $\bm{x}_R(\bm{\omega}_i,t)$ of \eqref{eq:dyn_xR} satisfies for all $\tau_{i-1} \leq t \leq \tau_{i-1}+T$ and all $\bm{\omega}_i \in \{1,2,\ldots, 2^M \}^K$
\begin{multline}
    |\bm{x}_R(\bm{\omega}_i,t)| \\ 
    \leq C_T e^{\mu_R(t-\tau_{i-1})} \left( |\bm{x}_{i-1}| + |\bm{u}_R(\bm{\omega}_i)|_{L^2(\tau_{i-1},t; \mathbb{R}^q)}\right). \label{eq:bound_muRxR}
\end{multline}
\end{lemma}

\begin{proof}
    Differentiate $|\bm{x}_R(\bm{\omega}_i,t)|^2$ using \eqref{eq:dyn_xR}, use \eqref{eq:def_muR}, integrate from $\tau_{i-1}$ to $t$, apply Cauchy-Schwarz in $L^2(\tau_{i-1},t; \mathbb{R}^q)$ and then Gronwall's lemma. 
\end{proof}

Our analysis will use Riccati theory. For the infinite-horizon OCP \eqref{eq:Jinf}--\eqref{eq:dyn_x}, let $P_\infty$ denote the (unique) symmetric positive-definite solution of the Algebraic Riccati Equation (ARE)
\begin{align}
    A^\top P_{\infty} + P_{\infty}A - P_{\infty}BW^{-1}B^\top P_{\infty} + Q = 0.
    \label{eq:ARE}
\end{align}
It is then well-known that, see, e.g., \cite[Section 5.1]{sage1968}, 
\begin{equation}
    \bm{u}^*_\infty(t) = -W^{-1} B^\top P_\infty \bm{x}_\infty^*(t). \label{eq:uinfPinf}
\end{equation}
Therefore, $\bm{x}^*_\infty(t)$ follows the dynamics generated by
\begin{equation}
    A_\infty = A-BW^{-1}B^\top P_\infty, \label{eq:Ainfty}
\end{equation}
which is stable, i.e. there exist $M_\infty \geq 1$ and $\mu_\infty > 0$
\begin{equation}
\| e^{A_\infty t} \| \leq M_{\infty} e^{-\mu_{\infty}t}. \label{eq:Ainfty_stab}
\end{equation}

For the finite-horizon OCP \eqref{eq:JT}--\eqref{eq:dyn_xT},  $t \in [0,T]$, let $P_T(t)$ solve the Riccati Differential Equation (RDE)
\begin{align}
-\dot{P}_T(t) & = A^\top P_T(t) + P_T(t)A \nonumber \\ & - P_T(t)BW^{-1}B^\top P_T(t) + Q, \quad P_T(T) = F, \label{eq:RDE}
\end{align}
on $t \in [0,T]$. It is well-known that, see, e.g., \cite[Section 5.2]{sage1968}, 
\begin{multline}
\bm{u}_T^*(t; \bm{x}_{i-1},\tau_{i-1}) = \\
-W^{-1}B^\top P_T(t - \tau_{i-1}) \bm{x}_T^*(t,\bm{x}_{i-1}, \tau_{i-1}). \label{eq:uTPT}
\end{multline}
For the randomized OCP \eqref{eq:JR}--\eqref{eq:dyn_xR}, let $P_R(\bm{\omega}_i,t)$ solves the Randomized Riccati Differential Equation (RRDE) on $[0,T]$
\begin{align} 
\dot{P}_R(\bm{\omega}_i, t) = A_R(\bm{\omega}_i, t)^\top P_R(\bm{\omega}_i, t) + P_R(\bm{\omega}_i, t)A_R(\bm{\omega}_i, t) \nonumber \\
-P_R(\bm{\omega}_i,t)BW^{-1}B^\top P_R(\bm{\omega}_i, t) + Q, \  P_R(\bm{\omega}_i, T) = F. \label{eq:RRDE}
\end{align}
Similarly as in \eqref{eq:uTPT}, it holds that
\begin{multline}
\bm{u}_R^*(\bm{\omega}_i,t; \bm{x}_{i-1},\tau_{i-1}) = \\
-W^{-1}B^\top P_R(\bm{\omega}_i,t - \tau_{i-1}) \bm{x}_R^*(\bm{\omega}_i,t,\bm{x}_{i-1}, \tau_{i-1}). \label{eq:uRPR}
\end{multline}


The following lemma shows that $P_T(t) \rightarrow P_\infty$ for $T \rightarrow \infty$. 
\begin{lemma} \label{lem:Pconv}
If $(A,B)$ is stabilizable, $(A,Q)$ is detectable, and $\mu_\infty$ is as in \eqref{eq:Ainfty_stab}, then for all $t \in [0,T]$
    \begin{equation}\label{P_infty - P_T 0}
        ||P_T(t) - P_{\infty}|| \leq C \| F - P_\infty \| e^{-2\mu_{\infty}(T-t)}.
    \end{equation}
\end{lemma}

\begin{proof}
See \cite{callier1994}. A shorter proof for the case that $(A,B)$ is controllable and $(A,Q)$ is observable is given in \cite{porretta2013}. 
\end{proof}

\begin{remark} \label{rem:PT}
Because $\| P_\infty \| \leq C$, Lemma \ref{lem:Pconv} implies that $\| P_T(t) \| \leq C$. 
\end{remark}


Let $V$ be a vector space. The expected value of a random variable $X : \{1,2,\ldots,2^M \}^K \rightarrow V$ depending on $\bm{\omega}_i$ is
\begin{equation}
     \mathbb{E}_i[X] = \hspace{-5mm} \sum_{\bm{\omega}_i \in \{1,2, \ldots, 2^M \}^K} \hspace{-5mm} X(\bm{\omega}_i) p(\bm{\omega}_i),
\end{equation}
where $p(\bm{\omega}_i) = p_{\omega_{i,1}} p_{\omega_{i,2}} \cdots p_{\omega_{i,K}}$. 
The expected value of a random variable $X : \{1,2,\ldots,2^M \}^{iK} \rightarrow V$ 
is denoted by
 \begin{multline}
     \mathbb{E}[X] = \hspace{-5mm}\sum_{\substack{\bm{\omega}_1, \bm{\omega}_2, \ldots , \bm{\omega}_i \\ \in \{1,2, \ldots, 2^M \}^K}}\hspace{-5mm} X(\bm{\omega}_1, \bm{\omega}_2, \ldots , \bm{\omega}_i) p(\bm{\omega}_1) p(\bm{\omega}_2) \cdots p(\bm{\omega}_i).
 \end{multline}
For the expected value of a random variable $X(\Omega_i)$ w.r.t.\ the last $\bm{\omega}_i$, we write $\mathbb{E}_i[X(\Omega_{i-1})]$ to indicate that the result depends on $\Omega_{i-1}$. 
For random variables $X(\Omega_i)$ and $Y(\Omega_i)$, 
\begin{align}
\mathbb{E}[XY] &\leq \sqrt{\mathbb{E}[X^2]\mathbb{E}[Y^2]}, \\
\sqrt{\mathbb{E}[(X+Y)^2]} &\leq \sqrt{\mathbb{E}[X^2]} + \sqrt{\mathbb{E}[Y^2]}.
\end{align}
Similar expressions hold for the expectation $\mathbb{E}_i$.

For $0 \leq s \leq t \leq T$, let $S_R(\bm{\omega}_i, t,s)$ be the evolution operator 
generated by $A_R(\bm{\omega}_i,t)$, i.e. $S_R(\bm{\omega}_i, t,s)\bm{x}= \bm{x}(\bm{\omega}_i, t)$ where 
\begin{equation}
\dot{\bm{x}}(\bm{\omega}_i,t) = A_R(\bm{\omega}_i, t) \bm{x}(\bm{\omega}_i,t), \qquad  \bm{x}(\bm{\omega}_i,s) = \bm{x}. 
\end{equation}

The following lemma from \cite{rbm} then shows that $S_R(\bm{\omega}_i,t,s)$ is (in expectation) close to $e^{A(t-s)}$ when $h \mathrm{Var}[A_R]$ is small. 
\begin{lemma} \label{lem:semigroup}
Let $\mathrm{Var}[A_R]$ and $\mu_R$ be as in \eqref{eq:def_varA} and \eqref{eq:def_muR} and $0 \leq s \leq t \leq T$, then
\begin{equation}
\mathbb{E}_i[\| S_R(t,s) - e^{A(t-s)} \|^2 ] \leq C_T e^{2 \mu_R (t-s)} h \mathrm{Var}[A_R].
\end{equation}
\end{lemma}
\begin{proof}
See \cite[Theorem 1 and Corollary 1]{rbm}. 
\end{proof}

With Lemma \ref{lem:semigroup}, it is also possible to bound the difference between controlled state trajectories. 

\begin{lemma} \label{lem:xRyR}
Let $\bm{u}_R: \{1,2, \ldots, 2^M \}^K \times [\tau_{i-1}, \tau_{i-1} + T] \rightarrow \mathbb{R}^q$ be a random control and $\bm{x}_{i-1} : \{1,2, \ldots, 2^M \}^K \rightarrow \mathbb{R}^n$ a random initial condition. If $\bm{x}_R(\bm{\omega}_i,t)$ and $\bm{y}_R(\bm{\omega}_i,t)$ satisfy
\begin{align}
\dot{\bm{x}}_R(\bm{\omega}_i, t) &= A\bm{x}_R(\bm{\omega}_i, t) + B \bm{u}_R(\bm{\omega}_i,t), \label{eq:lem_dynxR} \\
\dot{\bm{y}}_R(\bm{\omega}_i, t) &= A_R(\bm{\omega}_i,t-\tau_{i-1}) \bm{y}_R(\bm{\omega}_i, t) + B \bm{u}_R(\bm{\omega}_i,t), \label{eq:lem_dynyR} \\
\bm{x}_R(\bm{\omega}_i, &\tau_{i-1}) = \bm{y}_R(\bm{\omega}_i,\tau_{i-1}) = \bm{x}_{i-1}(\bm{\omega}_i), \label{eq:lem_xRyRinit}
\end{align}
then
\begin{align}
&\mathbb{E}[|\bm{y}_R(t;\bm{x}_i, \tau_i) - \bm{x}_R(t;\bm{x}_i, \tau_i)|^2] \leq C_T e^{2 \mu_R (t-\tau_i)} \times \\
&h \mathrm{Var}[A_R] \left(\max_{\bm{\omega}_i} |\bm{x}_{i-1}(\bm{\omega}_i)| + \max_{\bm{\omega}_i} |\bm{u}_R(\bm{\omega}_i)|_{L^2(\tau_{i-1},t;\mathbb{R}^q))} \right)^2. \nonumber 
\end{align}
\end{lemma}
\begin{proof}
By a slight modification of \cite[Theorem 2]{rbm} in which the random initial condition was not considered. 
\end{proof}

Although Riccati theory will be used in the analysis, the OCPs in Section \ref{sec:rbm_mpc} are typically more efficiently solved by a gradient-based algorithm, especially when $n$ is large. 
\section{Stability Analysis} \label{sec:stab}
For the stability result in Theorem \ref{thm:1} at the end of this section, we first establish two lemmas. Consider $i \in \mathbb{N}$ and $t \in [\tau_{i-1}, \tau_i]$. Because $\bm{x}_{R-M}(\Omega_{i}, t)$ satisfies \eqref{eq:dyn_yR},
\begin{equation}
\dot{\bm{x}}_{R-M}(\Omega_{i},t) = A_\infty\bm{x}_{R-M}(\Omega_{i},t) + \bm{r}(\Omega_i, t), \label{eq:dyn_xRM1}
\end{equation}
with $A_\infty$ as in \eqref{eq:Ainfty} and
\begin{multline}
\bm{r}(\Omega_i,t) = BW^{-1}B^\top P_\infty \bm{x}_{R-M}(\Omega_i,t) + B \bm{u}^*_R(\bm{\omega}_i,t) \\
= B W^{-1}B^\top (P_\infty - P_T(t-\tau_{i-1})) \bm{x}_{R-M}(\Omega_i,t) +  \\
B \underbrace{ \left(
W^{-1}B^\top P_T(t-\tau_{i-1}) \bm{x}_{R-M}(\Omega_i,t) + \bm{u}^*_R(\bm{\omega}_i,t) \right)}_{=:\bm{g}(\Omega_i,t)}, \label{eq:r_split}
\end{multline}
where $\bm{u}_R^{*}(\bm{\omega}_{i},t)$ denotes $\bm{u}_R^*(\bm{\omega}_i,t; \bm{x}_{R-M}(\Omega_{i-1}, \tau_{i-1}), \tau_{i-1})$ for brevity. 
The first auxiliary lemma is now as follows. 

\begin{lemma} \label{lem:PRPT}
Let $P_R(\bm{\omega}_i,t)$ and $P_T(t)$ satisfy \eqref{eq:RRDE} and \eqref{eq:RDE}, then, for $t \in [0,T]$,
\begin{equation}
\mathbb{E}_i[\|P_R(t) - P_T(t)\|] \leq C_T e^{2 \mu_R T}   \sqrt{h \mathrm{Var}[A_R]}. \label{eq:lem_PRPT}
\end{equation}
\end{lemma}
\begin{proof}
We will only proof \eqref{eq:lem_PRPT}  for  $t = 0$. The result for $t > 0$ can be obtained similarly. By definition,
\vspace{-1.11mm}
\begin{align}
&\| P_R(\bm{\omega}_i,0) - P_T(0) \| \\
&=  (\bar{\bm{x}}(\bm{\omega}_i))^\top (P_R(\bm{\omega}_i,0) - P_T(0)) \bar{\bm{x}}(\bm{\omega}_i) \nonumber \\
& = |J_R(\bm{\omega}_i, \bm{u}_R^*(\bm{\omega}_i); \bar{\bm{x}}(\bm{\omega}_i), 0) - 
J_T(\bm{u}^*_T(\bm{\omega}_i); \bar{\bm{x}}(\bm{\omega}_i),0)|, \nonumber
\end{align}
where $\bar{\bm{x}}(\bm{\omega}_i) = \mathrm{argmax}_{|\bm{x}| = 1} |\bm{x}^\top (P_R(\bm{\omega}_i,0) - P_T(0)) \bm{x}|$, $(\bm{u}_R^*(\bm{\omega}_i), \bm{x}_R^*(\bm{\omega}_i))$ and $(\bm{u}^*_T(\bm{\omega}_i), \bm{x}_T^*(\bm{\omega}_i))$ denote the control-state pairs that minimize $J_R(\bm{\omega}_i, \cdot\,; \bar{\bm{x}}(\bm{\omega}_i), 0)$ and $J_T(\cdot\,; \bar{\bm{x}}(\bm{\omega}_i), 0)$, respectively. We also write $J_R(\bm{\omega}_i, \cdot)$ and $J_T(\cdot)$ for $J_R(\bm{\omega}_i, \cdot\,; \bar{\bm{x}}(\bm{\omega}_i), 0)$ and $J_T(\cdot\,; \bar{\bm{x}}(\bm{\omega}_i), 0)$, respectively. 
Now introduce a random control $\bm{u}_{RT}(\bm{\omega}_i, t)$ by setting it equal to 
$\bm{u}_T^*(\bm{\omega}_i, t)$ when $J_T(\bm{u}^*_T(\bm{\omega}_i)) \leq J_R(\bm{\omega}_i, \bm{u}^*_R(\bm{\omega}_i))$ 
and equal to $\bm{u}_R^*(\bm{\omega}_i, t)$ otherwise and define $\bm{x}_{RT}(\bm{\omega}_i,t)$ and $\bm{y}_{RT}(\bm{\omega}_i,t)$ as the solutions of \eqref{eq:lem_dynxR}--\eqref{eq:lem_xRyRinit} on $[0, T]$ resulting from the control $\bm{u}_{RT}(\bm{\omega}_i, t)$ and the IC $\bar{\bm{x}}(\bm{\omega}_i)$. When $J_T(\bm{u}^*_T(\bm{\omega}_i)) \leq J_R(\bm{\omega}_i, \bm{u}^*_R(\bm{\omega}_i))$, 
\begin{align}
\| P_R (\bm{\omega}_i,0) - &P_T(0) \| = J_R(\bm{\omega}_i, \bm{u}^*_R(\bm{\omega}_i)) - J_T(\bm{u}^*_T(\bm{\omega}_i)) \nonumber \\ 
&\leq J_R(\bm{\omega}_i, \bm{u}_{RT}(\bm{\omega}_i)) - J_T(\bm{u}_{RT}(\bm{\omega}_i)), \label{eq:lem_PRPT_step1}
\end{align}
because $\bm{u}_{RT}(\bm{\omega}_i) = \bm{u}_T^*(\bm{\omega}_i)$ and because $\bm{u}^*_R(\bm{\omega}_i)$ minimizes $J_R(\bm{\omega}_i,\cdot)$. 
Similarly, when $J_R(\bm{\omega}_i, \bm{u}^*_R(\bm{\omega}_i)) < J_T(\bm{u}^*_T)$
\begin{align}
\| P_R (\bm{\omega}_i,0) - &P_T(0) \| = J_T(\bm{u}^*_T(\bm{\omega}_i)) - J_R(\bm{\omega}_i, \bm{u}^*_R(\bm{\omega}_i)) \nonumber \\ 
&\leq J_T(\bm{u}_{RT}(\bm{\omega}_i)) - J_R(\bm{\omega}_i, \bm{u}_{RT}(\bm{\omega}_i)), \label{eq:lem_PRPT_step3}
\end{align}
because $\bm{u}_{RT}(\bm{\omega}_i) = \bm{u}_R^*(\bm{\omega}_i)$ and $\bm{u}^*_T(\bm{\omega}_i)$ minimizes $J_T(\cdot)$. Combining \eqref{eq:lem_PRPT_step1} and \eqref{eq:lem_PRPT_step3} thus shows that
\begin{align}
    \| &P_R (\bm{\omega}_i,0) - P_T(0) \| \leq |J_R(\bm{\omega}_i, \bm{u}_{RT}(\bm{\omega}_i)) - J_T(\bm{u}_{RT}(\bm{\omega}_i))| \nonumber \\
&\leq
||\bm{y}_{RT}(\bm{\omega}_i,T)|_F^2 - |\bm{x}_{RT}(\bm{\omega}_i,T)|_F^2| \nonumber \\
&\quad + |\langle \bm{y}_{RT}(\bm{\omega}_i), Q \bm{y}_R(\bm{\omega}_i) \rangle_{L^2} - \langle \bm{x}_{RT}(\bm{\omega}_i), Q \bm{x}_{RT}(\bm{\omega}_i) \rangle_{L^2}| \nonumber \\
&\leq |\langle 2 \bm{x}_{RT}(\bm{\omega}_i,T) + \bm{e}_{RT}(\bm{\omega}_i,T), F \bm{e}_{RT}(\bm{\omega}_i, T) \rangle |  \nonumber \\
&\quad + |\langle 2 \bm{x}_{RT}(\bm{\omega}_i) + \bm{e}_{RT}(\bm{\omega}_i), Q \bm{e}_{RT}(\bm{\omega}_i) \rangle_{L^2}| \nonumber \\
&\leq 2 |\bm{x}_{RT}(\bm{\omega}_i,T)|_F |\bm{e}_{RT}(\bm{\omega}_i,T)|_F  + |\bm{e}_{RT}(\bm{\omega}_i,T)|_F^2 \nonumber  \\
&\quad + 2 \sqrt{\langle \bm{x}_{RT}(\bm{\omega}_i), Q \bm{x}_{RT}(\bm{\omega}_i) \rangle_{L^2} \langle \bm{e}_{RT}(\bm{\omega}_i), Q \bm{e}_{RT}(\bm{\omega}_i) \rangle_{L^2}} \nonumber \\
&\quad + \langle \bm{e}_{RT}(\bm{\omega}_i), Q \bm{e}_{RT}(\bm{\omega}_i) \rangle_{L^2},  
\label{eq:lem_PRPT_step6}
\end{align}
where $\langle \cdot, \cdot \rangle_{L^2}$ denotes the $L^2$-inner product on $[0, T]$ and 
$\bm{e}_{RT}(\bm{\omega}_i,t) = \bm{y}_{RT}(\bm{\omega}_i,t) - \bm{x}_{RT}(\bm{\omega}_i,t)$. Furthermore, when $J_T(\bm{u}^*_T(\bm{\omega}_i)) \leq J_R(\bm{\omega}_i, \bm{u}^*_R(\bm{\omega}_i))$,
\begin{align}
|\bm{x}_{RT}&(\bm{\omega}_i,T)|_F^2 + \langle \bm{x}_{RT}(\bm{\omega}_i), Q \bm{x}_{RT}(\bm{\omega}_i) \rangle_{L^2}  + \alpha |\bm{u}_{RT}(\bm{\omega}_i)|_{L^2}^2 \nonumber \\ 
&\leq J_T(\bm{u}^*_T(\bm{\omega}_i)) = |\bar{\bm{x}}(\bm{\omega}_i)|_{P_T(0)}^2 \leq C |\bar{\bm{x}}(\bm{\omega}_i))|^2,\label{eq:lem_PRPT_step2} 
\end{align}
where $\alpha > 0$ is the smallest eigenvalue of $W$ and it has been used that $\| P_T(t) \| \leq C$, see Remark \ref{rem:PT}. 
When $J_R(\bm{\omega}_i, \bm{u}^*_R(\bm{\omega}_i)) < J_T(\bm{u}^*_T)$
\begin{multline}
|\bm{x}_{RT}(\bm{\omega}_i,T)|_F^2 +
\langle \bm{x}_{RT}(\bm{\omega}_i), Q \bm{x}_{RT}(\bm{\omega}_i) \rangle_{L^2} + \alpha |\bm{u}_{RT}(\bm{\omega}_i)|^2_{L^2} \\ \leq J_R(\bm{\omega}_i, \bm{u}^*_R(\bm{\omega}_i)) 
\leq J_T(\bm{u}^*_T(\bm{\omega}_i)) \leq C |\bar{\bm{x}}(\omega_i)|^2.  \label{eq:lem_PRPT_step4}
\end{multline}
Because $|\bar{\bm{x}}(\bm{\omega}_i)| = 1$ and because \eqref{eq:lem_PRPT_step2} and \eqref{eq:lem_PRPT_step4} show that 
$|\bm{x}_{RT}(\bm{\omega},T)|_F \leq C$ and $\langle \bm{x}_{RT}(\bm{\omega}), Q \bm{x}_{RT}(\bm{\omega}) \rangle_{L^2} \leq C$, taking the expectation in \eqref{eq:lem_PRPT_step6} thus shows that
\begin{multline}
    \mathbb{E}_i[\| P_R(0) - P_T(0) \|] \leq C \mathbb{E}_i[|\bm{e}_{RT}(T)|] + C \mathbb{E}_i[|\bm{e}_{RT}(T)|^2] \\
    + C \mathbb{E}_i[|\bm{e}_{RT}|_{L^2}] + C \mathbb{E}_i[|\bm{e}_{RT}|^2_{L^2}].
    \label{eq:lem_PRPT_step7}
\end{multline}
Because $|\bar{\bm{x}}(\bm{\omega}_i)| = 1$ and because \eqref{eq:lem_PRPT_step2} and \eqref{eq:lem_PRPT_step4} show that 
$|\bm{u}_{RT}(\bm{\omega}_i)|_{L^2} \leq C$, applying Lemma \ref{lem:xRyR} shows that 
\begin{align}
\mathbb{E}_i[|\bm{e}_{RT}(t)|] &\leq C_T e^{2\mu_RT} h \mathrm{Var}[A_R], 
\end{align}
which also implies that $\mathbb{E}_i[ |\bm{e}_{RT} |_{L^2}^2] \leq C_T e^{2\mu_R T} h \mathrm{Var}[A_R]$. 
Because $\mathbb{E}_i[X] \leq \sqrt{\mathbb{E}_i[X^2]}$, \eqref{eq:lem_PRPT_step7} now shows that
\begin{multline}
    \mathbb{E}_i[\| P_R(0) - P_T(0) \|] \\
    \leq C_T e^{2 \mu_R T} \left(\sqrt{h \mathrm{Var}[A_R]}  + h \mathrm{Var}[A_R]\right).
\end{multline}
As we are interested in the limit $h\mathrm{Var}[A_R] \rightarrow 0$, we assume that $h \mathrm{Var}[A_R] \leq C$ (cf.\ the first paragraph of Section \ref{sec:preliminaries}) and \eqref{eq:lem_PRPT} follows because $h \mathrm{Var}[A_R] \leq C\sqrt{h \mathrm{Var}[A_R]}$. 
\end{proof}

With this result, $\mathbb{E}_i[| \bm{g}(\Omega_{i-1},t) |]$ can be bounded as follows. 

\begin{lemma} \label{lem:bound_g}
Let $\bm{g}(\Omega_i,t)$ be as in \eqref{eq:r_split}, then
\begin{multline}
\mathbb{E}_i[ |\bm{g}(\Omega_{i-1}, t)| ] \\ \leq
C_T e^{\mu_R (2T+\tau)} \sqrt{h \mathrm{Var}[A_R]} |\bm{x}_{R-M}(\Omega_{i-1}, \tau_{i-1})|. \label{eq:lem_bound_g}
\end{multline}
\end{lemma}

\begin{proof}
By using \eqref{eq:uRPR}, the expression for $\bm{g}(\Omega_i,t)$ in \eqref{eq:r_split} can be rewritten as
\begin{multline}
    \bm{g}(\Omega_i,t) = W^{-1}B^\top P_T(t-\tau_{i-1})(\bm{x}_{R-M}(\Omega_{i},t)-\bm{x}^*_R(\bm{\omega}_{i},t)) \\
    + W^{-1}B^\top (P_T(t-\tau_{i-1}) - P_R(\bm{\omega}_{i-1}, t-\tau_{i-1})) \bm{x}^*_R(\bm{\omega}_{i},t). \nonumber
\end{multline}
Taking the norm and expectation w.r.t.\ $\bm{\omega}_i$ yields
\begin{multline}
    \mathbb{E}_i[|\bm{g}(\Omega_{i-1},t)|] \leq C \mathbb{E}_i[|\bm{x}_{R-M}(\Omega_{i-1},t)-\bm{x}^*_R(t)|] \\
     + C \mathbb{E}_i[\| P_T(t-\tau_{i}) - P_R(t-\tau_{i})\|] \max_{\bm{\omega}_i} |\bm{x}^*_R(\bm{\omega}_{i},t)|. \label{eq:lem_g_step3}
\end{multline}
Lemma \ref{lem:PRPT} gives a bound for $\mathbb{E}_i[\| P_T(t-\tau_{i}) - P_R(t-\tau_{i})\|]$. 

For $\mathbb{E}_i[|\bm{x}_{R-M}(\Omega_{i-1},t)-\bm{x}^*_R(t)|]$, apply Lemma \ref{lem:xRyR} with $\bm{u}_R(\bm{\omega}_i,t) = \bm{u}_R^*(\bm{\omega}_i,t; \bm{x}_{R-M}(\Omega_{i-1}, \tau_{i}), \tau_{i})$ and the initial condition $\bm{x}_i(\bm{\omega}_i) = \bm{x}_{R-M}(\Omega_{i-1},\tau_i)$. This makes
\begin{align*}
\bm{y}_R(\bm{\omega}_i,t) &= \bm{x}^*_R(\bm{\omega}_i,t; \bm{x}_{R-M}(\Omega_{i-1}, \tau_{i}), \tau_{i}), \\
\bm{x}_R(\bm{\omega}_i,t) &= \bm{y}_R^*(\bm{\omega}_i,t; \bm{x}_{R-M}(\Omega_{i-1}, \tau_{i}), \tau_{i}) = \bm{x}_{R-M}(\Omega_i,t).
\end{align*}
The initial condition $\bm{x}_{R-M}(\Omega_{i-1},\tau_{i-1})$ does not depend on $\bm{\omega}_i$. Let $\alpha$ again denote the smallest eigenvalue of $W$, then
\begin{align}
\alpha &|\bm{u}^*_R(\bm{\omega}_i)|_{L^2}^2 \leq J_R(\bm{\omega}_i, \bm{u}_R^*(\bm{\omega}_i); \bm{x}_{R-M}(\Omega_{i-1}, \tau_{i}), \tau_i) \nonumber \\ 
&\leq J_R(\bm{\omega}_i, 0; \bm{x}_{R-M}(\Omega_{i-1}, \tau_{i}), \tau_i) \nonumber \\
&\leq C |\bm{x}_R(\bm{\omega}_i; \bm{x}_{R-M}(\Omega_{i-1}, \tau_{i}), \tau_i)|_{L^2(\tau_{i-1}, \tau_{i-1}+T; \mathbb{R}^q)}^2 \nonumber \\
&\leq C_Te^{2\mu_R T}|\bm{x}_{R-M}(\Omega_{i-1}, \tau_{i-1})|^2, \label{eq:boundJR}
\end{align}
where $\bm{x}_R(\bm{\omega}_i,t; \bm{x}_i, \tau_i)$ satisfies \eqref{eq:dyn_xR} with $\bm{u}_R(t) = 0$ and $\bm{x}_{i-1} = \bm{x}_{R-M}(\Omega_{i-1}, \tau_{i-1})$ and the last inequality follows from Lemma \ref{lem:bound_muRxR}. Lemma \ref{lem:xRyR} now shows that
\begin{align}
\mathbb{E}_i[| &\bm{x}_{R-M}(\Omega_{i-1},t)-\bm{x}^*_R(t)|] \nonumber \\
& \leq 
\sqrt{\mathbb{E}_i[|\bm{x}_{R-M}(\Omega_{i-1},t)-\bm{x}^*_R(t)|^2]} \nonumber \\ 
& \leq C_T  e^{\mu_R \tau} \sqrt{h \mathrm{Var}[A_R]} e^{\mu_R T} |\bm{x}_{R-M}(\Omega_{i-1},\tau_{i-1})|. \label{eq:bound_xRMxR}
\end{align}

To bound $|\bm{x}^*_R(\bm{\omega}_{i},t)|$, note that 
$\bm{x}^*_R(\bm{\omega}_{i},t)$ satisfies \eqref{eq:dyn_xR} with $\bm{u}_R(\bm{\omega}_i,t) = \bm{u}_R^*(\bm{\omega}_i,t)$ and $\bm{x}_{i-1} = \bm{x}_{R-M}(\Omega_{i-1}, \tau_{i-1})$. Inserting \eqref{eq:boundJR} into \eqref{eq:bound_muRxR} thus shows that for $t \in [\tau_{i-1}, \tau_i]$
\begin{equation}
|\bm{x}_R^*(\bm{\omega}_i,t)|
\leq C_T e^{\mu_R(T+\tau)} |\bm{x}_{R-M}(\Omega_{i-1},\tau_{i-1})|. \label{eq:bound_xR}
\end{equation}
Now insert \eqref{eq:bound_xRMxR}, \eqref{eq:lem_PRPT}, and \eqref{eq:bound_xR} into \eqref{eq:lem_g_step3} to find \eqref{eq:lem_bound_g}. 
\end{proof}

We are now ready to prove the main stability result. 

\begin{theorem}\label{thm:1}
If $(A,B)$ is stabilizable, $(A,Q)$ is detectable, and $M_\infty$ and $\mu_\infty$ are as in \eqref{eq:Ainfty_stab},  then
\begin{align}\label{eq:mr_stab}
    \mathbb{E}[|\bm{x}_{R-M}(t)|] \leq M_\infty e^{-\mu_{R-M} t} |\bm{x}_0|,
\end{align}
where
\begin{multline}
    \mu_{R-M} = \mu_\infty - C \| F - P_\infty \| e^{-2\mu_\infty(T-\tau)} \\
     - C_T e^{\mu_R(2T+\tau)}e^{\mu_\infty \tau} \sqrt{h \mathrm{Var}[A_R]} . \label{eq:def_muRM}
\end{multline}
\end{theorem}

\begin{proof}
Applying the variation of constants formula to \eqref{eq:dyn_xRM1}, taking the norm and the expectation yields
\begin{multline}
\mathbb{E}[|\bm{x}_{R-M}(t)|] \leq M_\infty e^{-\mu_\infty t} |\bm{x}_0| \\
+ M_\infty \int_0^t e^{-\mu_\infty (t-s)} \mathbb{E}[|\bm{r}(s)|] \ \mathrm{d}s, \label{eq:bound_xRM_step1}
\end{multline}
where \eqref{eq:Ainfty_stab} has been used. Taking the norm and the expectation (first w.r.t.\ $\bm{\omega}_{\lfloor s / \tau \rfloor+1}$ and then w.r.t.\ to the other $\bm{\omega}_j$'s) in \eqref{eq:r_split} using  Lemmas \ref{lem:Pconv} and \ref{lem:bound_g}, it follows that
\begin{equation}
\mathbb{E}[ |\bm{r}(s)| ] \leq C_1 \mathbb{E}[|\bm{x}_{R-M}(s)|] + C_2 \mathbb{E}[|\bm{x}_{R-M}(\tau_{\lfloor s/\tau \rfloor})|], \label{eq:r_bound}
\end{equation}
where we have introduced $C_1 = C \| F - P_\infty \| e^{-2\mu_\infty (T-\tau)}$ and $C_2 = C_T e^{\mu_R (2T + \tau)} \sqrt{h \mathrm{Var}[A_R]}$. By inserting \eqref{eq:r_bound} into \eqref{eq:bound_xRM_step1} and writing $f(t) = \mathbb{E}[\bm{x}_{R-M}(t)]$, we obtain
\begin{multline}
f(t) \leq M_\infty e^{-\mu_\infty t} |\bm{x}_0| \\
+ \int_0^t e^{-\mu_\infty (t-s)} \left( C_1 f(s) + C_2 f(\tau_{\lfloor s / \tau \rfloor})\right) \ \mathrm{d}s. 
\end{multline}
Setting $\hat{f}(t) = e^{\mu_\infty t} f(t)$, it follows that $\hat{f}(t) \leq \hat{F}(t)$ where
\begin{equation}
\hat{F}(t) = M_\infty |\bm{x}_0| + \int_0^t \left(C_1 \hat{f}(s) + C_2 e^{\mu_\infty \tau} \hat{f}(\tau_{\lfloor s / \tau \rfloor}) \right) \ \mathrm{d}s. \nonumber
\end{equation}
Because $\hat{F}(t)$ is monotonically increasing and $\hat{f}(t) \leq \hat{F}(t)$,
\begin{equation}
\hat{F}(t) \leq M_\infty |\bm{x}_0| + \left( C_1 + C_2 e^{\mu_\infty \tau} \right) \int_0^t \hat{F}(s) \ \mathrm{d}s. 
\end{equation}
By Gronwall's lemma, we thus obtain that
\begin{equation}
e^{\mu_\infty t} f(t) = \hat{f}(t) \leq \hat{F}(t) \leq M_\infty |\bm{x}_0| e^{(C_1 + C_2e^{\mu_\infty \tau})t}, 
\end{equation}
and the result follows. 
\end{proof}

\begin{remark}
Note that $\mu_{R-M} > 0$ will be positive for $h\mathrm{Var}[A_R]$ and $\|F - P_\infty \| e^{-2 \mu_\infty(T - \tau)}$ sufficiently small. 
\end{remark}

\begin{remark}
When $\mu_{R-M} > 0$,  the RBM-MPC strategy is stabilizing with probability 1. To see this, note that Markov's inequality 
and Theorem \ref{thm:1} imply that for any $\varepsilon > 0$
\begin{equation}
\mathbb{P}[|\bm{x}_{R-M}(t)| \geq \varepsilon] \leq \frac{\mathbb{E}[|\bm{x}_{R-M}(t)|]}{\varepsilon} \leq \frac{M_\infty e^{-t \mu_{R-M}} |\bm{x}_0| }{\varepsilon}. \nonumber
\end{equation}
Because $\mu_{R-M} > 0$, the probability that $\bm{x}_{R-M}(t)$ is outside any $\varepsilon$-neighborhood of the origin approaches zero for $t \rightarrow \infty$. 
\end{remark}



\section{Convergence} \label{sec:conv}
We first consider the convergence of MPC. Note that $\bm{x}_M(t)$ follows the dynamics generated by the $\tau$-periodic matrix
\begin{equation}
A_\tau(t) = A - BW^{-1}B^\top P_T(t \mod \tau). 
\end{equation}
We then have the following lemma. 
\begin{lemma} \label{lem:MPC}
If $(A,B)$ is stabilizable, $(A,Q)$ is detectable and $M_\infty$ and $\mu_\infty$ are as in \eqref{eq:Ainfty_stab}, then for all $0 \leq s \leq t$
\begin{equation}
\| e^{\int_s^t A_\tau(\sigma) \ \mathrm{d}\sigma} \| \leq M_\infty e^{-\mu_M t}, \label{eq:lem_MPC_1}
\end{equation}
where
\begin{equation}
\mu_M = \mu_\infty - C \| F - P_\infty \| e^{-2\mu_\infty (T-\tau)}. \label{eq:def_muM}
\end{equation}
Furthermore, if $\mu_M > 0$, then
\begin{multline}
|\bm{x}_M(t) - \bm{x}^*_\infty(t)| + 
|\bm{u}_M(t) - \bm{u}^*_\infty(t)| \\
 \leq C \| F - P_\infty\| e^{-2\mu_\infty(T-\tau)} |\bm{x}_0|. \label{eq:lem_MPC_2}
\end{multline}
\end{lemma}

\begin{remark}
Lemma \ref{lem:MPC} shows that the dynamics generated by $A_\tau(t)$ is stable for $T - \tau$ sufficiently large or $\| F - P_\infty\|$ sufficiently small and that $(\bm{x}_M(t), \bm{u}_M(t)) \rightarrow (\bm{x}^*_\infty(t), \bm{u}^*_\infty(t))$ for $T - \tau \rightarrow \infty$ or $\| F - P_\infty \| \rightarrow 0$. 
\end{remark}

\begin{proof} 
Let $\bm{x}(t)$ denote the solution to
 $\dot{\bm{x}}(t) = A_\tau(t) \bm{x}(t)$ with initial condition $\bm{x}(s) = \bm{x}_s$. By \eqref{eq:Ainfty}, 
\begin{equation}
\dot{\bm{x}}(t) = \left( A_\infty + B W^{-1} B^\top (P_\infty - P_T(t \mod \tau)) \right)\bm{x}(t). \label{eq:lem_MPC_3}
\end{equation}
The variation of constants formula thus shows that
\begin{multline}
\bm{x}(t) = e^{A_\infty(t-s)} \bm{x}_s \\
+ \int_s^t e^{A_\infty (t-\sigma)} BW^{-1}B^\top (P_\infty - P_T(\sigma \mod \tau)) \bm{x}(\sigma) \ \mathrm{d}\sigma. \nonumber
\end{multline}
Taking norms using \eqref{eq:Ainfty_stab} and Lemma \ref{lem:Pconv}, it follows that
\begin{multline}
|\bm{x}(t)| = M_\infty e^{-\mu_\infty(t-s)} |\bm{x}_s| \\
 + C \| F - P_\infty \| e^{-2\mu_\infty(T-\tau)} \int_s^t |\bm{x}(\sigma)| \ \mathrm{d}\sigma. 
\end{multline}
Applying Gronwall's lemma and noting that the initial condition $\bm{x}_s$ is arbitrary now yields \eqref{eq:lem_MPC_1}. For the bound on $\bm{e}_M(t) := \bm{x}_M(t) - \bm{x}^*_\infty(t)$ in \eqref{eq:lem_MPC_2}, note that $\dot{\bm{x}}_\infty^*(t) = A_\infty \bm{x}_\infty^*(t)$ and that $\bm{x}_M(t)$ satisfies \eqref{eq:lem_MPC_3}, so that
\begin{equation}
\dot{\bm{e}}_M(t) = A_\infty \bm{e}_M(t) - BW^{-1}B^\top (P_T(t \mod \tau) - P_\infty) \bm{x}_M(t). \nonumber
\end{equation}
Applying the variation of constants formula and taking the norm using \eqref{eq:Ainfty_stab}, Lemma \ref{lem:Pconv}, and the inequality $|\bm{x}_M(t)| \leq M_\infty |\bm{x}_0| \leq C |\bm{x}_0|$ when $\mu_M \geq 0$ by \eqref{eq:lem_MPC_1}, it follows that
\begin{equation}
|\bm{e}_M(t)| \leq C \| F - P_\infty \| e^{2\mu_\infty(T-\tau)} \int_0^t e^{\mu_\infty(t-s)} \ \mathrm{d}s |\bm{x}_0|. 
\end{equation}
The bound on $\bm{e}_M(t)$ follows because the remaining integral is bounded by $1/\mu_\infty \leq C$. For $\bm{u}_M(t) - \bm{u}_\infty(t)$, note that \eqref{eq:uTPT} implies that
\begin{equation}
    \bm{u}_M(t) = -W^{-1}B^\top P_T(t \mod \tau) \bm{x}_M(t), \label{eq:uMPT} 
\end{equation}
so that subtracting \eqref{eq:uinfPinf} shows that
\begin{multline}
\bm{u}_M(t) - \bm{u}_\infty(t) =
 W^{-1}B^\top (P_\infty - P_T(t \mod \tau))\bm{x}_M(t) \\
 -W^{-1}B^\top P_\infty \bm{e}_M(t).  
\end{multline}
Using Lemma \ref{lem:Pconv} and that $|\bm{x}_M(t)| \leq M_\infty |\bm{x}_0| \leq C |\bm{x}_0|$ for the first term, and the previously derived bound for $|\bm{e}_M(t)|$ for the second, the result follows. 
\end{proof}

Now the convergence of RBM-MPC can be established. 
\begin{theorem} \label{thm:2}
If $(A,B)$ is stabilizable, $(A,Q)$ is detectable, and $\mu_{R-M}$ in \eqref{eq:def_muRM} is positive, then
\begin{multline}
	\mathbb{E}[|\bm{x}_{R-M}(t) - \bm{x}_M(t)|] + \mathbb{E}[|\bm{u}_{R-M}(t) - \bm{u}_M(t)|]
	 \\
	 \leq \frac{C_T}{\mu_M} e^{\mu_R (2T + \tau)} \sqrt{h \mathrm{Var}[A_R]}|\bm{x}_0|.
\end{multline}
\end{theorem}
\begin{proof} 
Consider $i \in \mathbb{N}$ and $t \in [\tau_{i-1}, \tau_i)$.
For the bound on $\bm{e}_{R-M}(\Omega_i,t) = \bm{x}_{R-M}(\Omega_i,t) - \bm{x}_M(t)$, note that $\bm{x}_M(t)$ satisfies \eqref{eq:lem_MPC_3}, so that \eqref{eq:r_split} into \eqref{eq:dyn_xRM1} and subtracting \eqref{eq:lem_MPC_3} yields
\begin{equation}
\dot{\bm{e}}_{R-M}(\Omega_i,t) = A_\tau(t) \bm{e}_{R-M}(\Omega_i,t) + B\bm{g}(\Omega_i,t),
\end{equation}
and $\bm{e}_{R-M}(\Omega_i,0) = 0$. Applying the variation of constants formula, taking the norm and the expectation thus shows that
\begin{equation}
\mathbb{E}[|\bm{e}_{R-M}(t)| ] = C \int_0^t \left\| e^{\int_s^t A_\tau(\sigma) \ \mathrm{d}\sigma} \right\| \mathbb{E}[|\bm{g}(s)|] \ \mathrm{d}s. \label{eq:bound_eRM_step0}
\end{equation}
By Lemma \ref{lem:bound_g}, it follows that
\begin{align}
\mathbb{E}[|\bm{g}(s)|] &\leq C_T e^{\mu_R(2T + \tau)} \sqrt{h \mathrm{Var}[A_R]} \mathbb{E}[|\bm{x}_{R-M}(\tau_{\lfloor s/\tau \rfloor})|] \nonumber \\
&\leq C_T e^{\mu_R(2T + \tau)} \sqrt{h \mathrm{Var}[A_R]} |\bm{x}_0| \label{eq:bound_expectation_g2}
\end{align}
where it has been used that $\mathbb{E}[|\bm{x}_{R-M}(t)|] \leq M_\infty |\bm{x}_0| \leq C |\bm{x}_0|$ by Theorem \ref{thm:1} because $\mu_{R-M} \geq 0$. 
Using \eqref{eq:lem_MPC_1} and \eqref{eq:bound_expectation_g2} in \eqref{eq:bound_eRM_step0}, the bound for $\bm{e}_{R-M}(\Omega_i, t)$ follows because the integral of $e^{-\mu_{M}(t-s)}$ is bounded by $1/\mu_M$. 

To bound $\bm{u}_{R-M}(\Omega_i, t) - \bm{u}_M(t)$, note that for $t \in [\tau_{i-1}, \tau_i)$
\begin{align}
 \bm{u}_{R-M}(\Omega_i, t) &= \bm{u}_R^*(\bm{\omega}_i, t; \bm{x}_{R-M}(\Omega_{i-1}, \tau_{i-1}), \tau_{i-1}).
\end{align}
Subtracting \eqref{eq:uMPT} using the definition of $\bm{g}(\Omega_i,t)$ in \eqref{eq:r_split} yields
\begin{multline}
\bm{u}_{R-M}(\Omega_i, t) - \bm{u}_M(t) = \\
\bm{g}(\Omega_i,t) - W^{-1}B^\top P_T(t \mod \tau) \bm{e}_{R-M}(\Omega_i,t). 
\end{multline}
The bound now follows after taking the norm and the expected value, and then using Lemma \ref{lem:bound_g} to bound $\mathbb{E}[|\bm{g}(t)|]$ and the previously derived estimate for $\mathbb{E}[|\bm{e}_{R-M}(t)|]$. 
\end{proof}

\begin{remark}
Combining Theorem \ref{thm:2} and \eqref{eq:lem_MPC_2}, one obtains estimates for $\mathbb{E}[|\bm{x}_{R-M}(t) - \bm{x}^*_\infty(t)|] + \mathbb{E}[|\bm{u}_{R-M}(t) - \bm{u}^*_\infty(t)|]$. 
\end{remark}

The estimates also indicate a natural approach to tuning the parameters in RBM-MPC. First, $T-\tau$ should be chosen such that the MPC strategy is stabilizing with sufficient margin, i.e. such that $C \| F - P_\infty\| e^{2\mu_\infty(T-\tau)} \ll \mu_\infty$. After that, $h$ can be chosen such that $\mu_{R-M} > 0$ and such that RBM-MPC leads to a sufficiently good approximation of MPC.

\section{Numerical Example} \label{sec:numerical}
We consider a problem of the form \eqref{eq:Jinf}--\eqref{eq:dyn_x} with $n \in \{ 11, 101, 1001\}$ states and $m = 1$ input, $A \in \mathbb{R}^{n \times n}$ is
\begin{equation}
A = (n-1)^2 \begin{bmatrix}
-2 & 1 & 0 & \cdots  & 0 & 1 \\
1 & -2 & 1 &  & 0 & 0 \\
0 & 1 & -2 &  & 0 & 0 \\
\vdots & & & \ddots & \\
0 & 0 & 0 &  & -2 & 1 \\
1 & 0 & 0 &  & 1 & -2 
\end{bmatrix}, \label{eq:def_Aexample}
\end{equation}
the first $(n-1)/10$ entries of $B \in \mathbb{R}^{n \times 1}$ are 1 and the others are zero, $Q = I_n/(n-1)$, and $R = 1$, where $I_n$ denotes the $n \times n$ identity matrix. The matrix $A$ in \eqref{eq:def_Aexample} could for example be obtained as a spatial discretization of a heat equation on a circle.
The time interval $[0, \infty)$ is truncated to $[0,200]$ and discretized by the Crank-Nicholson scheme with a step size $\Delta t = 1$. Note that each time step is at least of $O(n^2)$, because the matrix $I_n - \tfrac{\Delta t}{2}A$ is not tridiagonal. The OCPs are solved by a steepest descent algorithm, in which the gradients are computed using the adjoint state, see \cite{Apel}, and the stepsize minimizes the functional in the direction of the gradient. The algorithm is stopped when the relative change in the control is below $10^{-5}$ or after 1000 iterations. 


To construct the randomized matrix $A_R(\bm{\omega}_i,t)$, note that $A$ can be written as the sum of $M = n$ interconnection matrices as in \eqref{eq:Asplit}, where the first $n-1$ interconnection matrices $A_m$ are zero except for a diagonal block of the form
\begin{equation}
(n-1)^2 \begin{bmatrix}
-1 & 1 \\ 1 & -1
\end{bmatrix}, 
\end{equation}
and the last interconnection matrix has only nonzero entries in its four corners. The sum of the first $n-1$ interconnection matrices leads to a tridiagonal matrix, which reduces the computational cost for each time step to $O(n)$, see, e.g., \cite[Section 2.1.1]{quarteroni1994}.
In fact, the symmetry of the problem implies that omitting \emph{any} one of the $n$ submatrices $A_m$ reduces the computational cost for one time step to $O(n)$. 
A probability $1/n$ is assigned to each subset of $\{1,2, \ldots, n \}$ of size $n-1$. The probabilities $\pi_m$ in \eqref{eq:pim} are thus $\pi_m = \frac{n-1}{n}$. The grid spacing $h$ is chosen as small as possible, so $h = \Delta t$. 
All $A_m$ are dissipative, so $\mu_R = 0$ by Remark \ref{rem:dissipative}.

%

\begin{figure*}
\centering
\subfloat[The controls $u_{R-M}(\Omega_i,t)$, $u_M(t)$, and $u^*_\infty(t)$. \label{fig:results_time_control}]{%
\includegraphics[width=0.45\textwidth]{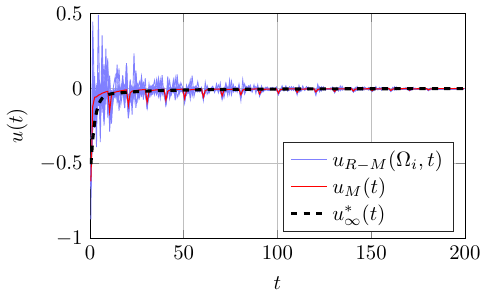}}
\hfill
\subfloat[The norm of the state trajectories $\bm{x}_{R-M}(\Omega_i,t)$, $\bm{x}_M(t)$, and $\bm{x}^*_\infty(t)$. \label{fig:results_time_state}]{%
\includegraphics[width=0.45\textwidth]{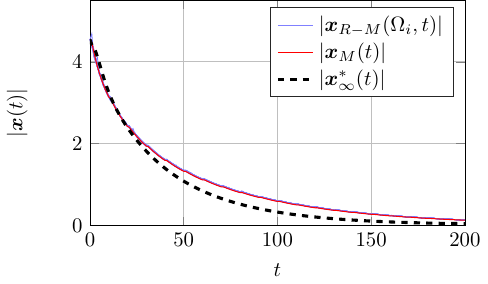}}
\caption{The RBM-MPC control and state trajectory $u_{R-M}(\Omega_i,t)$ and $\bm{x}_{R-M}(\Omega_i, t)$ for 20 realizations of $\Omega_i$ compared to $u_M(t)$, $\bm{x}_M(t)$, $u^*_\infty(t)$, and $\bm{x}^*_\infty(t)$ for $n=100$, $h=1$, $\tau = 10$, and $T = 15$. The lines for $|\bm{x}_{R-M}(\Omega_i,t)|$ and $|\bm{x}_M(t)|$ in Figure \ref{fig:results_time_state} almost overlap. 
}
\label{fig:results_time}
\end{figure*}

Figure \ref{fig:results_time_control} compares 20 realizations of the RBM-MPC control $u_{R-M}(\Omega_i,t)$ to the MPC control $u_M(t)$ and the infinite horizon control $u^*_\infty(t)$ for $n = 100$ spatial grid points. As can be seen, $u^*_\infty(t)$ is smooth, $u_M(t)$ jumps when $t$ is a multiple of $\tau = 10$, and the realizations of $u_{R-M}(\Omega_i,t)$ contain high-frequent oscillations related to the grid spacing $\Delta t = h = 1$. 
Figure \ref{fig:results_time_state} shows that that despite the relatively large deviations of $u_{R-M}(\Omega_i,t) $ from $u_M(t)$, $|\bm{x}_{R-M}(\Omega_i,t)|$ is very close to $|\bm{x}_M(t)|$ for all 20 considered realizations $\Omega_i$. RBM-MPC thus leads to almost the same decay rate as the MPC here. Note that $T = 15$ is not much larger than $\tau = 10$, but the simulations indicate that MPC and RBM-MPC are stabilizing.

Table \ref{tab:times_varyn} shows that the running times for RBM-MPC are smaller than those for MPC, which are again smaller than those for solving the OCP on $[0,200]$ directly. The numbers between round brackets in Table \ref{tab:times_varyn} indicate the estimated standard deviation of the running times based on 20 runs. For $n = 100$, MPC is almost 3 times faster than a classical optimal control approach, and RBM-MPC is again almost 3 times faster than MPC. For $n = 1000$, MPC is still approximately 3 times faster than solving the OCP directly, but RBM-MPC is 5 times faster than MPC. Note that the relative speed-up of RBM-MPC compared to MPC may not always match theoretical estimates due to overhead and potential additional iterations in the RBM-constrained OCP compared to the original OCP.

These observations are particularly interesting because Table \ref{tab:error_varyn} shows that the errors do not increase significantly when $n$ is increased. The numbers between round brackets in Table \ref{tab:error_varyn} indicate the estimated standard deviation based on 20 realizations of $\Omega_i$. Here, $\| \bm{x} \|_{L^\infty} := \max_t \sqrt{(\bm{x}(t))^\top \bm{x}(t)/n}$. 

\begin{table}
\caption{Running times for a varying number of spatial grid points $n$ \\ ($h = 1$, $T = 15$, $\tau = 10$)}
\centering
\begin{tabular}{|l||r|r|r|} \hline
Running times [s]  & \multicolumn{1}{|c|}{$n=10$} & \multicolumn{1}{|c|}{$n=100$} & \multicolumn{1}{|c|}{$n=1000$} \\ \hline \hline
Optimal Control    & 12.9 ($\pm$1.35) & 32.7 ($\pm$1.53) & 218.5 ($\pm$4.36) \\ \hline
MPC          & 4.6 ($\pm$0.44) & 11.0 ($\pm$0.66) & 70.5 ($\pm$4.51) \\ \hline
RBM-MPC      & 2.0 ($\pm$0.23)  &  3.6 ($\pm$0.71) & 14.1 ($\pm$1.33)  \\ \hline
\end{tabular}
\label{tab:times_varyn}
\end{table}

\begin{table}
\caption{Errors for a varying number of spatial grid points $n$ \\ ($h = 1$, $T = 15$, $\tau = 10$)}
\centering
\begin{tabular}{|l||c|c|c|} \hline
Relative errors [-]  & $n=10$ & $n=100$ & $n=1000$ \\ \hline \hline
$| u_{R-M} - u^*_\infty |_{L^2}$ & 0.76 ($\pm$0.28) & 0.59 ($\pm$0.23) & 0.53 ($\pm$0.11) \\ \hline
$| u_{R-M} - u_M |_{L^2}$ & 0.63 ($\pm$0.30) & 0.41 ($\pm$0.27) & 0.33 ($\pm$0.16) \\ \hline
$| u_M   - u^*_\infty |_{L^2}$ & 0.41 ($\pm$0.00) & 0.37 ($\pm$0.00) & 0.39 ($\pm$0.00) \\ \hline \hline
$\| \bm{x}_{R-M} - \bm{x}^*_\infty \|_{L^\infty}$ & 0.35 ($\pm$0.09) & 0.28 ($\pm$0.08) & 0.25 ($\pm$0.04) \\ \hline
$\| \bm{x}_{R-M} - \bm{x}_M \|_{L^\infty}$ & 0.17 ($\pm$0.08) & 0.11 ($\pm$0.07) & 0.08 ($\pm$0.04)  \\ \hline
$\| \bm{x}_M   - \bm{x}^*_\infty \|_{L^\infty}$ & 0.24 ($\pm$0.00) & 0.22 ($\pm$0.00) & 0.22 ($\pm$0.00) \\ \hline
\end{tabular}
\label{tab:error_varyn}
\end{table}

\begin{figure*}
\centering
\subfloat[Varying $h$, $T = 15$, $\tau = 10$. \label{fig:state_conv_h}]{%
\includegraphics[width=0.28\textwidth]{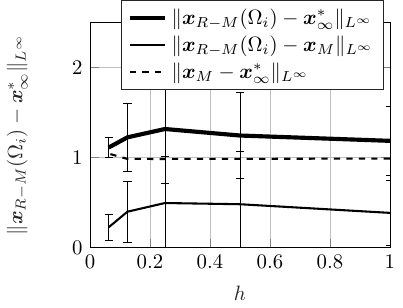}}
\hfill
\subfloat[$h = 1$, varying $T$, $\tau = 10$. \label{fig:state_conv_T}]{%
\includegraphics[width=0.28\textwidth]{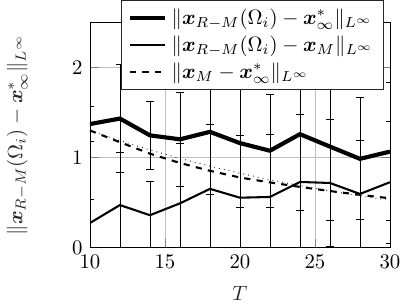}}
\hfill
\subfloat[$h = 1$, $T = 40$, varying $\tau$. \label{fig:state_conv_tau}]{%
\includegraphics[width=0.28\textwidth]{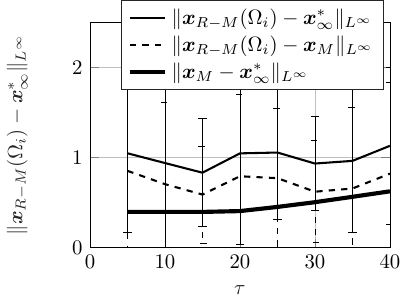}}
\caption{Differences between the RBM-MPC state trajectory $\bm{x}_{R-M}(\Omega_i, t)$, the MPC state trajectory $\bm{x}_M(t)$, and the infinite horizon state trajectory $\bm{x}^*_\infty(t)$ for $n =100$. The error bars indicate the $2\sigma$ confidence intervals estimated based on 20 realizations of $\Omega_i$. }
\label{fig:state_conv}
\end{figure*}

\begin{figure*}
\centering
\subfloat[Varying $h$, $T = 15$, $\tau = 10$. \label{fig:control_conv_h}]{%
\includegraphics[width=0.28\textwidth]{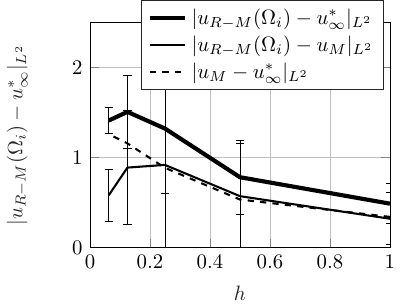}}
\hfill
\subfloat[$h = 1$, varying $T$, $\tau = 10$. \label{fig:control_conv_T}]{%
\includegraphics[width=0.28\textwidth]{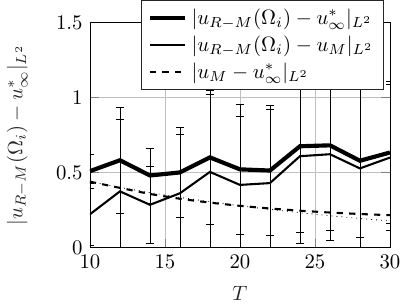}}
\hfill
\subfloat[$h = 1$, $T = 40$, varying $\tau$. \label{fig:control_conv_tau}]{%
\includegraphics[width=0.28\textwidth]{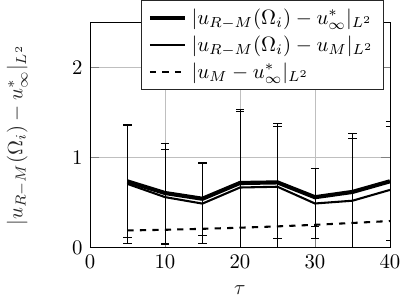}}
\caption{Differences between the RBM-MPC control $u_{R-M}(\Omega_i,t)$, the MPC control $u_M(t)$, and the infinite horizon control $u^*_\infty(t)$ for $n =100$. The error bars indicate the $2\sigma$ confidence intervals estimated based on 20 realizations of $\Omega_i$. }
\label{fig:control_conv}
\end{figure*}

The convergence rates from Lemma \ref{lem:MPC}  and Theorem \ref{thm:2} are validated in Figure \ref{fig:state_conv} and \ref{fig:control_conv}. Figures \ref{fig:state_conv_h} and \ref{fig:control_conv_h} show that $\|\bm{x}_{R-M}(\Omega_i) - \bm{x}_M \|_{L^\infty}$ and $|u_{R-M}(\Omega_i) - u_M|_{L^2}$ decay as $\sqrt{h}$ for $h \rightarrow 0$ and that $\bm{x}_{R-M}(\Omega_i)$ and $u_{R-M}(\Omega_i)$ do not converge to $\bm{x}^*_\infty$ and  $u^*_\infty$ for $h \rightarrow 0$, as the estimates from Section \ref{sec:conv} indicate. 
Figures \ref{fig:state_conv_T} and \ref{fig:control_conv_T} show that $\| \bm{x}_M - \bm{x}^*_\infty \|_{L^\infty}$ and $|u_M - u^*_\infty|_{L^2}$ are proportional to $e^{-2\mu_\infty T}$, as Lemma \ref{lem:MPC} indicates. Increasing $T$ increases $\| \bm{x}_{R-M}(\Omega_i) - \bm{x}_M \|_{L^\infty}$ and $|u_{R-M}(\Omega_i) - u_M|_{L^2}$, which confirms that the constant $C_T$ in Theorem \ref{thm:2} increases with $T$. 
Figures \ref{fig:state_conv_tau} and \ref{fig:control_conv_tau} show that varying $\tau$ does not affect $\|\bm{x}_{R-M}(\Omega_i) - \bm{x}_M \|_{L^\infty}$ and $|u_{R-M}(\Omega_i) - u_M|_{L^2}$ strongly and  $\|\bm{x}_{M} - \bm{x}^*_\infty \|_{L^\infty}$ and $|u_M - u^*_\infty|_{L^2}$ increase with $\tau$. 

The code used to generated the results in this section can be found on \url{https://github.com/DCN-FAU-AvH}. 

\section{Conclusion and Perspectives} \label{sec:conclusions}
This paper considers a randomized MPC strategy called RBM-MPC to efficiently approximate the solution of an large-scale infinite-horizon linear-quadratic OCP. In RBM-MPC, the finite-horizon OCPs in each MPC-iteration are simplified by replacing the system matrix $A$ by a randomized one. The estimates in this paper demonstrate that 1) RBM-MPC is stabilizing for $h \mathrm{Var}[A_R]$ sufficiently small and either $T - \tau$ sufficiently large or $\| F - P_\infty \|$ sufficiently small, and 2) RBM-MPC states and controls converge in expectation to their MPC counterparts for $h \mathrm{Var}[A_R] \rightarrow 0$. 
In an example with $n = 100$ states, RBM-MPC is 9 times faster than solving the OCP direcly and 3 times faster than classical MPC. 

The estimates in this note form a natural starting point for the analysis of RBM-MPC in nonlinear and/or constrained settings in future works. The computational advantage of RBM-MPC has already been demonstrated in a nonlinear setting, see \cite{ko2021model}. Because the training of residual Deep Neural Networks (DNNs) can be seen viewed as a nonlinear OCP (see, e.g.,  \cite{benning2019, esteve2020}), RBM-MPC may also be applied to speed up the training of DNNs. 
RBM-MPC may also be used for the control of (networks of) PDEs, that for example appear in the modeling of gas transport, see, e.g., \cite{herty2007}. 

Finally, other variations of RBM-MPC could be considered. One variation would be to first fix a RBM approximation over the whole time axis $[0, \infty)$ and use this as the plant model for MPC. 
Another interesting variation would be to consider a new (independent) RBM approximation in each step of the gradient descent algorithm used to solve the (finite horizon) OCPs in MPC. 

\bibliographystyle{IEEEtran}
\bibliography{bilbio}

\end{document}